\documentclass[eqno,11pt]{article}

\usepackage{amssymb}
\usepackage{euscript}

\usepackage{times}

\usepackage[latin1]{inputenc}
\usepackage{amsmath}
\usepackage{amsfonts}

 \evensidemargin\oddsidemargin
\begin{document}

\renewcommand{\theequation}{\thesection.\arabic{equation}}
\newtheorem{theorem}{Theorem}[section]
\newtheorem{lemma}[theorem]{Lemma}
\newtheorem{proposition}[theorem]{Proposition}
\newtheorem{corollary}[theorem]{Corollary}
\newtheorem{remark}[theorem]{Remark}
\newtheorem{fact}[theorem]{Fact}
\newtheorem{problem}[theorem]{Problem}
\newtheorem{example}[theorem]{Example}

\newcommand{\eqnsection}{
\renewcommand{\theequation}{\thesection.\arabic{equation}}
    \makeatletter
    \csname  @addtoreset\endcsname{equation}{section}
    \makeatother}
\eqnsection

\def\r{{\mathbb R}}
\def\e{{\mathbb E}}
\def\p{{\mathbb P}}
\def\q{{\mathbb Q}}
\def\qa{{\mathbb Q}^{(\alpha)}}
\def\z{{\mathbb Z}}
\def\n{{\mathbb N}}
\def\T{{\mathbb T}}
\def\L{{\mathcal  L}}
\def\S{{\mathcal S}}
\def\F{{\mathcal F}}
\def\G{{\mathcal G}}
\def\M{{\mathbb M}}

\vglue20pt

\centerline{ \Large \bf  The almost sure limits of the minimal position  and the additive  martingale  }

\medskip

\centerline{\Large  \bf   in a branching random walk }

\bigskip
\bigskip

\medskip

 \centerline{Yueyun Hu\footnote{Département de Mathématiques, Université Paris 13, Sorbonne Paris Cité,  LAGA   (CNRS UMR 7539),  F--93430 Villetaneuse.   Research partially  supported by ANR  2010 BLAN 0125 Email: yueyun@math.univ-paris13.fr}}

\medskip

\centerline{\it Universit\'e Paris XIII}

\bigskip

{\leftskip=2truecm
\rightskip=2truecm
\baselineskip=15pt
\small

\noindent{\slshape\bfseries Summary.}   Consider a   real-valued branching random walk in the boundary case.  Using   the techniques developed by   A{\"{\i}}d{\'e}kon and Shi \cite{AS12},  we give two integral tests which describe respectively the lower limits for the minimal position and the upper limits for the associated additive martingale. 
\medskip



} 

\bigskip
\bigskip

\section{Introduction}

Let  $\{ V(u), u \in \T\}$ be   a discrete-time  branching random walk on the real line $\r$, where $\T$ is an  Ulam-Harris tree which describes the genealogy of the particles and $V(u)\in \r$ is the position of the particle $u$. When  a particle $u$ is at $n$-th generation, we write $\vert u \vert=n$ for   $n \ge0$. The branching random walk $V$ can be described as follows:  At the beginning, there is a single particle $\varnothing$ located at $0$. The particle $\varnothing$ is  also the root of $\T$.  At the generation $1$, the root dies and  gives birth to some  point process $\L$ on $\r$. The point process $\L$ constitutes the first generation of the branching random walk $\{ V(u), \vert u \vert=1\}$. The next generations are defined  by recurrence:  For each $\vert u \vert=n$ (if such $u$ exists),   the  particle $u$  dies  at the  $(n+1)$-th generation   and  gives birth to    an independent copy of $\L$ shifted by $V(u)$.   The collection of all children of all $u$ together with their positions  gives the $(n+1)$-th generation.  The whole system may survive   forever  or die out after some generations.

Plainly   $\L= \sum_{ \vert u \vert =1} \delta_{ \{ V(u)\}}$.   Assume  $\e [\L(\r)]  >1$ and that  \begin{equation}\label{hyp}    \e \left [ \int e^{-x} \L(dx)\right]=1 , \qquad \e \left[ \int x e^{-   x} \L(dx)\right]=0. \end{equation}

When the hypothesis \eqref{hyp} is fulfilled, the branching random walk is   called  in the boundary case in the literature (see e.g. Biggins and Kyprianou \cite{BK04} and  \cite{BK05}, A{\"{\i}}d{\'e}kon and Shi \cite{AS12}).  Under some integrability  conditions, a general branching random walk can  be reduced to the boundary case after a linear transformation, see Jaffuel \cite{J12} for detailed discussions.  We shall assume \eqref{hyp} throughout this paper.

Denote by $\M_n:= \min_{ \vert u \vert = n} V(u)$  the minimal position of the branching random walk at generation $n$    (with convention $\inf\emptyset  \equiv \infty$).  Hammersly \cite{H74}, Kingman \cite{K75} and Biggins \cite{B76} established the law of large numbers for $\M_n$ (for any general branching random walk), whereas the second order limits have recently attracted many attentions, see \cite{AbR, HS09, BZ06, A11} and the references therein.  In particular,  A{\"{\i}}d{\'e}kon \cite{A11}     proved the convergence in law of $\M_n- {3\over 2} \log n$ under \eqref{hyp} and some mild conditions,  which gives a discrete    analog of Bramson \cite{B78}'s theorem on    the branching brownian motion.  

Concerning   the almost sure limits of $\M_n$,  there is  a phenomena of fluctuation at the logarithmic scale (\cite{HS09}): Under \eqref{hyp} and some extra  integrability assumption: $\exists \,\delta>0$ such that  $\e[ \L(\r)^{1+\delta}] < \infty$ and $\e \big[ \int_{\r} (e^{ \delta x } + e^{- (1+\delta)x } ) \L(dx)\big] < \infty$,  the following almost sure limits hold:  \begin{eqnarray*} \limsup_{ n \to \infty} { \M_n \over \log n} &=& { 3\over 2}, \qquad \p^*\mbox{-a.s.}, 
\\ \liminf_{ n \to \infty} { \M_n \over \log n} &=& { 1\over 2}, \qquad \p^*\mbox{-a.s.}, 
\end{eqnarray*}

\noindent where here and in the sequel,      $$ \p^*(\cdot):= \p \left( \cdot \vert \S \right),$$ and $\S $ denotes the event that the whole system survives.  The upper bound ${ 3\over 2} \log n$ is  the   usual  fluctuation for $\M_n$ because $\M_n - {3\over 2} \log n$ converges in law (\cite{A11}). It is  a natural question to ask how $\M_n$ can approach the unusual lower bound ${1\over2} \log n  $.

  A{\"{\i}}d{\'e}kon and Shi \cite{AS12}   proved that under \eqref{hyp} and the following integrability conditions \begin{eqnarray} \label{int1}
   && \sigma^2:=   \e \Big[ \int_{ \r} x^2 e^{-x} \L(dx)  \Big] < \infty, \\
   && \e \Big[ \eta ( \log_+ \eta)^2 + \widetilde \eta \log_+\widetilde \eta  \Big] < \infty, \label{int2}
\end{eqnarray} where   $\eta:= \int_{\r} e^{-x} \L(dx) $, $ \widetilde \eta:= \int_{0}^\infty  x\, e^{-x}    \L(dx) $   and $\log_+ x:= \max( 0, \log x)$,   then $$ \liminf_{ n\to \infty} \left( \M_n - { 1\over 2} \log n \right) = - \infty, \qquad \mbox{$\p^*$-a.s.}$$
 
 \noindent 
 Furthermore, they asked whether there is some deterministic sequence $a_n \to \infty$ such that $$ - \infty < \liminf_{ n \to \infty} { 1\over a_n}  \left( \M_n - { 1\over 2 }\log n \right) < 0, \qquad \mbox{$\p^*$-a.s.}?$$

 The answer is yes: we can choose $a_n= \log \log n$. Moreover, we can give an integral test to describe the lower limits of $\M_n$:

 \begin{theorem}\label{t:min}  Assume  \eqref{hyp}, \eqref{int1} and \eqref{int2}.  For any function $f\uparrow \infty$,   \begin{equation}\label{test1}  \p^* \left( \M_n - {1\over 2} \log n < -  f(n) , \quad \mbox{ i.o.} \right) =        \left\{
     \begin{array}{ll}
         0           \\ \\
    1
     \end{array} \right.  \,\, \Longleftrightarrow \,\,  \int^\infty { dt \over t \exp(f(t))}  \left\{
     \begin{array}{ll}
         < \infty          \\ \\
      = \infty  
     \end{array} \right. , \end{equation}  where $\mbox{i.o.}$ means infinitely often as the relevant index  $n \to \infty$. 
     \end{theorem}

 As a consequence  of the integral test \eqref{test1}, we have that for any $\varepsilon>0$, $\p^*$-a.s. for all large $n\ge n_0(\omega)$, $ \M_n - { 1 \over 2} \log n \ge - (1+\varepsilon) \log \log n $ whereas there exists infinitely often $n $ such that $ \M_n - { 1 \over 2} \log n \le -   \log \log n . $  Hence    $\p^*$-a.s., $\liminf_{n\to\infty} {1\over \log \log n}( \M_n - { 1 \over2} \log n ) = -1$.

 The behaviors of the minimal position $\M_n$ are  closely related to the so-called additive martingale $(W_n)_{n \ge0}$:  $$ W_n:= \sum_{ \vert u \vert =n} e^{ -V(u)}, \qquad n \ge0,$$

\noindent with the usual convention: $\sum_\emptyset \equiv0$. By Biggins \cite{B76} and Lyons \cite{L97}, $W_n \to 0$ almost surely as $n \to \infty$.  The problem to find the rate of convergence (or a Seneta-Heyde norming)  for $W_n$ arose  in  Biggins and Kyprianou \cite{BK04} and was studied in \cite{HS09}. A{\"{\i}}d{\'e}kon and Shi \cite{AS12}  gave a definite result  to this problem. Let  \begin{equation}\label{defdn} D_n:= \sum_{ \vert u \vert =n} V(u) e^{- V(u)}, \qquad n\ge 1, \end{equation} be the derivative martingale (which is a martingale under the boundary condition \eqref{hyp}). It was shown in Biggins and Kyprianou \cite{BK04} that $\p$-a.s., $D_n $ converges to some nonnegative random variable $D_\infty$. Moreover under \eqref{hyp}, \eqref{int1} and \eqref{int2}, $\p^*$-a.s., $D_\infty >0$, as shown in \cite{BK04} and \cite{A11}. 

\medskip
{\noindent \bf Theorem (A{\"{\i}}d{\'e}kon and Shi \cite{AS12}).}
 {\it   Assume \eqref{hyp}, \eqref{int1} and \eqref{int2}. Then under $\p^*$, $$ \sqrt{n} W_n \, \,  {\, {\stackrel{(p)}{\rightarrow}}\, }  \, \, \sqrt{ 2\over \pi \sigma^2}\,  D_\infty, $$ as $n \to \infty$.  Moreover $$ \limsup_{n \to \infty} \sqrt{n} \, W_n = \infty, \qquad\mbox{$\p^*$-a.s.}$$ }
 
Furthermore A{\"{\i}}d{\'e}kon and Shi   conjectured that \begin{equation}\label{conj} \liminf_{n \to \infty} \sqrt{n} \, W_n = \sqrt{ 2\over \pi \sigma^2}\, D_\infty, \qquad \mbox{$\p^*$-a.s.}\end{equation}

    The upper limits of $W_n$ can be described as follows:

\begin{theorem} \label{t:upper}  Assume \eqref{hyp}, \eqref{int1} and \eqref{int2}.  For any  function $f \uparrow \infty$, $\p^*$-almost surely,      \begin{equation}\label{test2} \limsup_{ n \to \infty} { {\sqrt n} \, W_n \over f(n)} =        \left\{
     \begin{array}{ll}
         0           \\ \\
      \infty  
     \end{array} \right.  \,\, \Longleftrightarrow \,\,  \int^\infty { dt \over t f(t)}  \left\{
     \begin{array}{ll}
         < \infty          \\ \\
      = \infty  
     \end{array} \right.  .\end{equation}
     \end{theorem}

Concerning the lower limits of $ W_n$, we confirm  \eqref{conj} under a stronger integrability assumption:   There exists some small constant  $\varepsilon_0>0$ such that  \begin{equation} \e \left[ \eta^{1+\varepsilon_0} + \int e^{-x} \vert x \vert^{2+\varepsilon_0} \L(dx) \right] < \infty. \label{int3}
\end{equation}

\noindent It is easy to see that the condition \eqref{int3} is stronger than  \eqref{int1} and \eqref{int2}.

\begin{proposition}  \label{p:lower}  Assume \eqref{hyp} and  \eqref{int3}.  We have $$ \liminf_{n \to \infty} \sqrt{n} \, W_n = \sqrt{ 2\over \pi \sigma^2}\,  D_\infty, \qquad \p^*\mbox{-a.s.}$$
\end{proposition}

Combining Theorems \ref{t:min} and \ref{t:upper}, we can roughly say that the main contribution to the upper limits of $    W_n$ comes from the term $e^{-\M_n}$.  According to Madaule \cite{M10},  and A{\"{\i}}d{\'e}kon, Berestycki, Brunet and Shi \cite{ABBS}, Arguin, Bovier and Kistler  \cite{ABK} (for the branching brownian motion), the branching random walk seen from the minimal position converges in law to some point process, in particular, $W_n e^{ \M_n}$ converges in law as $n \to \infty$, but   we are not able to determine  the almost sure fluctuations  of $W_n e^{ \M_n}$.

The  whole paper uses essentially  the techniques developed by  A{\"{\i}}d{\'e}kon and Shi \cite{AS12}.  
To show   Theorems \ref{t:min} and \ref{t:upper}, we firstly remark  that  both two theorems share the same integral test and that  since  $W_n \ge e^{- \M_n}$,   it is enough to prove the convergence part in the integral test \eqref{test2} and the divergence part in \eqref{test1}. The convergence part in \eqref{test2} will follow from an application  of Doob's maximal inequality to a certain martingale. To prove   the divergence part in \eqref{test1}, we  shall use  the arguments in A{\"{\i}}d{\'e}kon and Shi \cite{AS12} (the proof of their Lemma 6.3) to estimate  a second moment, then apply  Borel-Cantelli's lemma.   We can also directly prove Theorem \ref{t:upper}   without the use of the divergence part of \eqref{test1}. Finally, the proof of Proposition \ref{p:lower} relies on a result (Lemma \ref{l:AS41}) which is also implicitly contained  in  A{\"{\i}}d{\'e}kon and Shi \cite{AS12} (by following  the proof of their Proposition 4.1).

The rest of this paper is organized as follows: In Section \ref{s:2}, we recall some known results on the branching random walk (many-to-one formula, change of measure) and on a real-valued random walk. In Section \ref{s:thm2}, we prove Theorems  \ref{t:min} and \ref{t:upper}, whereas the proof of Proposition \ref{p:lower} will be given in Section \ref{s:4}.
 
 Throughout this paper, $f(n) \sim g(n)$ as $n \to \infty$ means that $ \lim_{n \to \infty} {f (n) \over g(n)} =1$ and $(c_i, 1\le i \le 36)$ denote some positive  constants. 
\section{Preliminaries} \label{s:2}

\subsection{Many-to-one formula for the branching random walk}

In this subsection, we recall some change  of measure formulas  in  the branching random walk, for the details we refer to \cite{BK04,CRW,LPP, AS12, Shi} and the references therein.  

At first let us fix some notations which will be used throughout this paper: For $ \vert u \vert = n$, we write $[\varnothing, u]\equiv \{ u_0:= \varnothing, u_1, ..., u_{n-1}, u_n =u\}$ the shortest path     from the root $\varnothing$ to $u$ such that $\vert u_i\vert =i$ for any $0\le i \le n$. For any $u, v \in \T$, we use the partial order  $ u < v$ if $u$ is an ancestor of $v$ and  $u \le v$ if $u < v$ or $u=v$.  We also denote by   ${\buildrel \leftarrow \over v}$   the parent of $v$.   

Under \eqref{hyp}, there exists a centered real-valued random walk $\{S_n, n\ge0\}$ such that for any $n \ge1$ and any measurable function  $f:\r^n \to \r_+$, \begin{equation} \label{many} \e \Big[ \sum_{ \vert u \vert =n} e^{- V(u)}  f( V(u_1), ..., V(u_n)) \Big]= \e \left[ f(S_1, ..., S_n)\right]. \end{equation} 

\noindent Moreover under \eqref{int1}, $  \sigma^2 = \mbox{Var}(S_1)= \e \Big[ \sum_{\vert u \vert =1} (V(u))^2 e^{- V(u)} \Big] \in (0, \infty).$ 

The renewal function $R(x)$ related to the random walk $S$  is defined as follows: \begin{equation}\label{rx}  R(x):= \sum_{k=0}^\infty \p \left( S_k \ge -x , \, S_k < \min_{ 0 \le j \le k-1} S_j \right) , \qquad x \ge0, \end{equation} and $R(x)=0$ if $x <0$.    Moreover,  \begin{equation}\label{defcr} \lim_{x \to \infty} { R(x)\over x}= c_R , \end{equation} with some positive constant $c_R$ (see Feller \cite{feller}, pp.612).  

For   $\alpha \ge0$, we define as in A{\"{\i}}d{\'e}kon and Shi  \cite{AS12} two truncated processes:  For any $n \ge0$, \begin{eqnarray}  W_n^{(\alpha)} &:= &\sum_{ \vert u \vert =n} e^{- V(u)} 1_{( \underline V(u) \ge -\alpha)}, \qquad \label{wna} \\
	D_n^{(\alpha)} &:=& \sum_{ \vert u \vert =n} R_\alpha(V(u)) e^{- V(u)} 1_{( \underline V(u) \ge -\alpha))}, \label{dna}
\end{eqnarray} where  $\underline V(u) := \min_{ \varnothing \le v \le u } V(v) $,  $R_\alpha(x):= R(\alpha+x)$ and $R$ is the renewal function defined in \eqref{rx}.

Denote by $(\F_n, n\ge0)$ the natural filtration of the branching random walk.  If the branching random walk starts from $V(\varnothing)=x$, then we denote its law by $\p_x$ (with $\p=\p_0$). According to Biggins and Kyprianou \cite{BK04}, $(D_n^{(\alpha)}, n\ge0)$ is a $(\p_x, (\F_n))$-martingale and on some enlarged probability space (more precisely on the space of marked trees enlarged by an infinite ray $(\xi_n, n \ge0)$, called spine), we may construct a family of probabilities $(\q_x^{(\alpha)}, x \ge -\alpha)$ such that for any $x\ge -\alpha$,  the following statements (i), (ii) and (iii) hold:

(i)  For all $n\ge1$,  
\begin{eqnarray}
   && { d \qa_x \over d \p_x } \big\vert _{\F_n}= { D_n^{(\alpha)}\over D_0^{(\alpha)}}, \label{bk1} \\
   && \q_x^{(\alpha)} \left( \xi_n = u \big\vert \F_n\right)= {1\over D_n^{(\alpha)}}\, R_\alpha(V(u)) e^{-V(u)} 1_{(\underline V(u) \ge -\alpha)}, \qquad \forall \vert u \vert =n. \label{bk2} \end{eqnarray}

(ii)  Under $\qa_x$, the process $\{V(\xi_n), n \ge0\}$ along the spine $(\xi_n)_{n\ge0}$,  is distributed as the random walk $(S_n, n\ge0)$ under $\p$ conditioned to stay in $[-\alpha, \infty)$.  Moreover for any $n \ge1, x\ge -\alpha$ and $f : \r^n \to \r_+$,     \begin{equation} \label{bk4} 
	  \e_{\qa_x}  \Big[  f(V(\xi_1), ..., V(\xi_n))\Big]		 = {1\over R_\alpha(x)}  \e_x  \Big[  f(S_1, ..., S_n)  R_\alpha(S_n) 1_{(\underline S_n \ge -\alpha)} \Big]. 
	\end{equation}

(iii) Let $\G_n:= \sigma\{ u, V(u): {\buildrel \leftarrow \over u} \in \{\xi_k, 0\le k < n\}\},  $ $n\ge0$. Under $\qa_x$ and conditioned on $\G_\infty$, for all $u \not \in \{ \xi_k, k\ge0\}$ but ${\buildrel \leftarrow \over u} \in \{\xi_k,  k \ge0\}$ the induced branching random walk $(V(uv), \vert v \vert \ge0)$ are independent and are distributed as $\p_{V(u)}$, where $\{uv, \vert v \vert \ge0\}$ denotes the subtree of $\T$ rooted at $u$.

\medskip

Let us mention   that as a consequence of (i), the following many-to-one formula holds:  For any $n \ge1, x\ge -\alpha$ and $f : \r^n \to \r_+$,  \begin{equation}
 \label{bk3} 
	\e_x \Big[  \sum_{ \vert u \vert =n} e^{- V(u)} R_\alpha(V(u))  f(V(u_1), ..., V(u_n)) 1_{ ( \underline V(u) \ge -\alpha)} \Big] =   R_\alpha(x) e^{-x}\,  \e_{\qa_x} \Big[ f(V(\xi_1), ..., V(\xi_n))\Big] .
\end{equation}

\subsection{Estimates on a centered real-valued random walk}

We collect here some estimates on a   real-valued random walk $\{S_k, k\ge0\}$, centered and with finite   variance $\sigma^2 >0$.  Let $\underline S_n:=  \min_{  0\le i \le n} S_i$, $\forall\, n \ge 0$. Recall \eqref{rx}  for the renewal function $R(\cdot)$.

\begin{fact}  There exists some constant $c_1>0$ such that for any    $x \ge  0, $  \begin{eqnarray}\label{F1} \p_x \Big(  \underline S_n  \ge 0 \Big)  &\le  & c_1\, (1+x) n^{-1/2}, \qquad \forall  \, n \ge1, \\
	\p_x \big( \underline S_{n-1} > S_n \ge0\big) &\le&  c_1 (1+x) R(x) n^{-3/2}, \qquad \forall  \, n \ge1,  \label{AJ} \\
	\label{K1} \p_x \Big( \underline S_n \ge 0\Big) &\sim & \theta\,  R(x) n^{-1/2}, \qquad \mbox{as } n \to \infty,  \end{eqnarray}  
	
	\noindent with $\theta= { 1 \over c_R} \sqrt{ 2 \over \pi \sigma^2}$.  Moreover there is   $c_2>0$ such that for any $b \ge a \ge 0, x \ge 0, n\ge 1$,   \begin{equation} \label{AS1} \p_x\Big( S_n \in [a, b],  \underline S_n  \ge 0\Big) \le c_2  (1+x) (1+b- a) (1+b) n^{-3/2},\end{equation}
For any $0< r< 1$, there exists some $c_3=c_3(r) >0$ such that for all  $b \ge a \ge 0, x, y \ge 0, n\ge 1$, \begin{equation}\label{AS2}  \p_x \left( S_n \in [y+a, y+b], \, \underline S_n  \ge 0, \min_{ r   n \le j \le n} S_j \ge y \right) \le c_3 \, (1+x) (1+b-a) (1+b) n^{-3/2}. \end{equation}
\end{fact}

See Feller (\cite{feller}, Theorem 1a, pp.415) for \eqref{F1},   A{\"{\i}}d{\'e}kon and Jaffuel (\cite{AJ11}, equation (2.8))  for \eqref{AJ}, A{\"{\i}}d{\'e}kon and Shi \cite{AS12} for \eqref{AS1} and  \eqref{AS2},  and  Kozlov \cite{K76} and  Lemma 2.1  in \cite{AS12} for \eqref{K1} with the identification of the constant $\theta= { 1 \over c_R} \sqrt{ 2 \over \pi \sigma^2}$.

We end this section by an estimate on the stability on $x$ in the convergence \eqref{K1}.

\begin{lemma}\label{l:rw}  Let $S$ be a centered random walk with positive variance. There exists a constant $c_4>0$ such that for all $n \ge1$ and $x \ge 0 $, $$ \Big\vert { \p_x( \underline S_n \ge 0)  \over  R(x) \p ( \underline S_n \ge0) } - 1\Big\vert \le c_4  { 1+x \over \sqrt n}.$$  \end{lemma}

{\noindent \bf Proof of Lemma \ref{l:rw}.}  Denote in this proof by $\varrho(n):= \p( \underline S_n \ge0)$ for $n\ge 0$.  
Let $x\ge 0$. By considering the first $k \in [0, n]$ such that $S_k= \underline S_n$, we get   that \begin{eqnarray*}
\p_x(\underline S_n \ge0) &=& \p_x(\underline S_n \ge x)+ \sum_{k=1}^n \p_x \Big( \underline S_{k-1} > S_k \ge0, \min_{ k < j \le n} S_j  \ge S_k\Big) \\
	&=& \varrho(n)+   \sum_{k=1}^n \p_x \Big( \underline S_{k-1} > S_k \ge0\Big) \varrho(n-k) ,
\end{eqnarray*}

\noindent by the Markov property at $k$.  Note that $R(x)=  1+\sum_{k=1}^\infty \p_x \big( \underline S_{k-1} > S_k \ge0\big)$. It follows that \begin{equation}  \p_x(\underline S_n \ge0)  \le  R(x) \varrho(n) +  \sum_{k=1}^n \p_x \Big( \underline S_{k-1} > S_k \ge0\Big) [ \varrho(n-k)- \varrho(n)],  \label{16a} \end{equation}

\noindent and \begin{equation}
	 \p_x(\underline S_n \ge0)   \ge   R(x) \varrho(n) -  \sum_{k=n+1}^\infty \p_x \Big( \underline S_{k-1} > S_k \ge0\Big)  \varrho(n) .  \label{16b} \end{equation}

Denote respectively by $I_{\eqref{16a}}$ and   $I_{\eqref{16b}}$ the sum $\sum_{k=1}^n$ in \eqref{16a} and the sum $\sum_{k=n+1}^\infty $ in \eqref{16b}.   Let  $T^-:= \inf\{ j \ge1: S_j <0\}$. By the local limit theorem    (Eppel \cite{E79}, see also \cite{VW09}, equation (22)), if the distribution of $S_1$ is non-lattice, then  \begin{equation} \label{eppel}
 \p \Big( T^- =k\Big)  \, \sim\, {C_- \over k^{3/2}}, \qquad k \to \infty, \end{equation}

 \noindent with some positive constant $C_-$. Moreover Eppel \cite{E79} mentioned that  a modification of \eqref{eppel} holds in the lattice distribution case. Then there exists some constant $c_5>0$ such that for all $k \ge1$,  \begin{equation} \label{eppel2}  \p \Big( T^- =k\Big) \le {c_5 \over k^{3/2}}. \end{equation}

\noindent 	It follows that for any $k \le n$, $ \varrho(n-k) - \varrho(n)= \p \big( n-k < T^- \le n\Big) \le c_5 \sum_{i=n-k+1}^n i^{-3/2}$. Then by \eqref{AJ},  $$
I_{\eqref{16a}}  \le  c_6 (1+x) R(x) \sum_{k=1}^n  k^{-3/2} \sum_{i=n-k+1}^n i^{-3/2}.   
$$

Elementary computations show that   $\sum_{k=1}^{n/2}  k^{-3/2} \sum_{i=n-k+1}^n i^{-3/2}   \le \sum_{k=1}^{n/2}  k^{-3/2} \times k ({n\over2})^{-3/2} =O({1\over n})$    and    $\sum_{k=n/2}^{n}  k^{-3/2} \sum_{i=n-k+1}^n i^{-3/2} \le ({n\over2})^{-3/2} \sum_{i=1}^n i^{-3/2} \times i = O({1\over n})$. Hence $ I_{\eqref{16a}} \le c_7 (1+x) R(x) {1\over  n}  \le c_8 (1+x) R(x) {1\over \sqrt{n}} \varrho(n) $ by   \eqref{K1}. 

Finally again by \eqref{AJ}, we get that $I_{\eqref{16b}} \le c_9 (1+x) R(x) { 1\over \sqrt n} \varrho(n)$.   Then the Lemma follows from \eqref{16a} and \eqref{16b}. $\Box$

\section{Proofs  of Theorems \ref{t:min} and  \ref{t:upper}}\label{s:thm2}

In view of the inequality: $W_n \ge e^{ - \M_n}$, the convergence part of the integral test   \eqref{test2} yields  that of \eqref{test1}, whereas the divergence part of the integral test \eqref{test1} implies that of \eqref{test2}. We only need to show the convergence  part in \eqref{test2} and the divergence part in \eqref{test1}.

\subsection{Proof of the convergence  part in Theorem \ref{t:upper}:}

\begin{lemma}\label{l:1}  Assume \eqref{hyp}. For any $\alpha \ge0$, there exists some constant  $c_{10}=c_{10}(\alpha)>0$  such that for any $1< n  \le m$ and $\lambda>0$, we have $$ \p \Big( \max_{ n  \le   k \le m}  \sqrt{k} W_k^{(\alpha)} > \lambda\Big) \le  c_{10}{ \log n \over \sqrt n} + c_{10} { 1 \over \lambda} \, \sqrt{ m\over n}.$$
\end{lemma}

{ \noindent\bf Proof of Lemma \ref{l:1}.}  For $n \le k \le m+1$, define $$  \widetilde W_k^{(\alpha, n)}:= \sum_{ \vert u \vert = k} e^{-V(u)} 1_{( \underline V( u_n) \ge -\alpha)}, $$

\noindent  where as before $\underline V(u_n):= \min_{ 1 \le j \le n} V(u_j)$ and $u_n$ is the ancestor of $u$ at $n$-th generation.  Then $\widetilde W_n^{(\alpha, n)}= W_n^{(\alpha)}$.  

For $k \in [n, m]$, $ \widetilde W^{(\alpha, n)}_{ k+1}= \sum_{ \vert v \vert =k}   1_{( \underline V(v_n) \ge -\alpha)} \sum_{  u: {\buildrel \leftarrow \over u} =v } e^{-  V(u) }.$  The branching property implies that $ \e \Big(\widetilde W_{k+1}^{(\alpha, n)}  \vert \F_k\Big)= \widetilde W_k^{(\alpha, n)}$ for   $k \in [n, m]$.  By Doob's maximal inequality, $$ \p\Big( \max_{ n \le k\le m} \sqrt k \widetilde W_k^{(\alpha, n)} \ge \lambda \Big) \le {  \sqrt{m} \over \lambda} \e (\widetilde W_m^{(\alpha, n)}) = { \sqrt m \over \lambda} \e ( W_n^{(\alpha)}).$$

By the   many-to-one formula  \eqref{many} and the random walk estimate \eqref{F1}, $$ \e (W_n^{(\alpha)}) = \p \Big( \underline S_n \ge -\alpha\Big) \le {c_{11} \over \sqrt n},$$ with $c_{11}:= c_1 (1+\alpha) $.   It follows that $$ \p\Big( \max_{ n \le k\le m} \sqrt k \widetilde W_k^{(\alpha, n)} \ge \lambda \Big) \le  {c_{11}   \over \lambda} \sqrt{ m \over n}.$$

\noindent Comparing $\widetilde W_k^{(\alpha, n)}$ and $  W_k^{(\alpha)}$, we get that  $$\p \Big( \max_{ n  \le   k \le m}  \sqrt{k} W_k^{(\alpha)} > \lambda\Big) \le \p\Big( \min_{ n \le k \le m} \min_{ \vert u \vert =k} V(u) < -\alpha\Big) + {c_{11} \over \lambda} \sqrt{ m \over n}.$$

\noindent The proof of the Lemma will be finished if we can show that for all $n \ge2$, \begin{equation}\label{31ab} \p\Big(   \min_{ \vert u \vert  \ge n } V(u) < -\alpha\Big)  \le c_{10} { \log n \over \sqrt n}. \end{equation}

\noindent To this end, let us apply the following    known result  (see e.g. \cite{Shi}):       $$ \p \Big( \inf_{ u \in \T} V(u) < -x \Big) \le   e^{-x}, \qquad 
\forall\, x\ge 0.$$ 

\noindent Then for all $n\ge 2$,  \begin{eqnarray*} && \p\Big( \min_{ k \ge n } \min_{ \vert u \vert =k} V(u) < -\alpha\Big)  
	\\&\le& \p \Big( \inf_{ u \in \T} V(u) < -   \log n \Big) + \p \Big(\min_{ k \ge n } \min_{ \vert u \vert =k} V(u) < -\alpha, \inf_{ v \in \T} V(v) \ge -   \log n\Big) \\
	& \le&    {1\over n}  + \sum_{ k=n}^\infty \e\Big[ \sum_{ \vert u \vert =k} 1_{( V(u) < -\alpha, V(u_n)\ge -\alpha, ..., V(u_{k-1}) \ge -\alpha,\,  \underline V(u) \ge -   \log n)} \Big] \\
	&=&   {1\over n}  + \sum_{ k=n}^\infty \e \Big[ e^{S_k} 1_{( S_k < -\alpha, S_n \ge -\alpha, ..., S_{k-1} \ge -\alpha, \, \underline S_k \ge -   \log n)}\Big]\\
	&\le&   {1\over n} + e^{-\alpha} \p\Big( \underline S_n \ge -   \log n\Big),
\end{eqnarray*} where the above equality is due to the many-to-one formula  \eqref{many}.  Using   \eqref{F1} to bound the above probability term, we get \eqref{31ab} and  the Lemma. $\Box$

 \medskip
{\noindent\bf Proof of the convergence part in Theorem \ref{t:upper}:} Let $f$ be nondecreasing such that $\int^\infty  { dt \over t f(t)} < \infty$. 
 Let $n_j:= 2^j$ for large $j \ge j_0$. Then $ \sum_{j=j_0}^\infty  {1\over f(n_j)} < \infty$.  By using Lemma \ref{l:1},  $$ \p\Big( \max_{ n_j \le k \le n_{ j+1}} \sqrt k W_k^{(\alpha)} > f(n_j) \Big) \le  c_{10} { \log n_j \over \sqrt{ n_j}} + c_{10}\,   \ { \sqrt{2}  \over f(n_j)}, $$ whose sum on $j$ converges.  The Borel-Cantelli lemma implies that $\p$-a.s.   for all large $k  $,  $\sqrt{k} W_k^{(\alpha)} \le f(k)$.  Replacing $f(k)$ by $\varepsilon f(k)$ with an arbitrary constant $\varepsilon>0$, we get that $$ \limsup_{ k \to \infty} { \sqrt k W_k^{(\alpha)} \over f(k)} = 0, \qquad \mbox{$\p$-a.s.},$$ 
 
 \noindent for any $\alpha\ge 0$. By considering a countable $\alpha\to \infty$  (for instance $\alpha$ integer) and by using the fact that  $ W^{(\alpha)}_k =  W_k$ on the set $\{ \inf_{u \in \T} V(u) \ge -\alpha\}$, we get the convergence  part. $\Box$.

 \subsection{Proof of the divergence  part in Theorem \ref{t:min}:}


 The following lemma is a slight modification of  A{\"{\i}}d{\'e}kon and Shi \cite{AS12}'s  Lemma 6.3: 
 
 \begin{lemma}[\cite{AS12}] \label{l:as} There exist some constants $K>0$ and $c_{12}=c_{12}(K)>0$ such that for all $n \ge 2, 0 \le \lambda \le {1 \over 3} \log n$,  \begin{equation}\label{efk1} c_{12}\, e^{- \lambda} \le  \p \Big( \bigcup_{ k=n+1}^{2n} \big( E^{(n, \lambda)}_k \cap F^{(n, \lambda)}_k \big) \not = \emptyset\Big)  \le {1\over c_{12}} e^{-\lambda}, \end{equation} where for $n < k\le 2n$, \begin{eqnarray*}
E^{(n, \lambda)}_k&:=& \Big\{ u: \vert u \vert=k, {1 \over  2} \log n - \lambda \le V(u) \le {1 \over  2} \log n - \lambda+ K, \, V(u_i) \ge a^{(n, \lambda)}_i, \forall \, 0 \le i \le k\Big\},	 \\
F^{(n, \lambda)}_k&:=& \Big\{ u: \vert u \vert=k, \sum_{ v \in \Upsilon(u_{i+1})} ( 1+ (V(v)- a_i^{(n,\lambda)})_+) e^{ - (V(v)- a_i^{(n,\lambda)})} \le K \, e^{- b_i^{(k, n)}}, \, \forall \, 0 \le i < k\Big\},  
\end{eqnarray*}   where for $u \in \T \backslash \{ \varnothing\}$,     $\Upsilon(u):= \{ v: v \not=u,  {\buildrel \leftarrow \over v}= {\buildrel \leftarrow \over u}\}$ denotes  the set of brothers of $u$, $x_+:=\max(x, 0)$,     $$ a_i^{(n, \lambda)}:=  \Big( {1 \over 2 }\log n - \lambda \Big)\, 1_{( { n \over 2} < i \le 2n)}, \qquad 0 \le i \le 2n, $$ and for $n < k \le 2n$, $$ b_i^{(k, n)}:= i^{ 1/12} 1_{( 0 \le i \le { n\over2})}+ (k-i)^{1/12} \, 1_{( {n\over2} < i \le k)}, \qquad 0 \le i \le k.$$ 
 \end{lemma}
 
 {\noindent\bf Proof of Lemma \ref{l:as}.}   The proof of the lower bound in  \eqref{efk1}  [by the second moment method]  goes in the same way as that of Lemma 6.3 in A{\"{\i}}d{\'e}kon and Shi \cite{AS12} [We also keep their        notations], by replacing ${ 1 \over 2} \log n$ in their proof by ${ 1 \over 2} \log n-\lambda$.  Moreover, a similar computation of the second moment will be given  in the proof of Lemma \ref{l:bc}. Then we omit  the details.
 
 The upper bound in \eqref{efk1} is a simple consequence of the many-to-one formula:   Defining  $s:= { 1 \over 2} \log n - \lambda$, we have that  \begin{eqnarray*} 
  \p \Big( \bigcup_{ k=n+1}^{2n}   E^{(n, \lambda)}_k  \not=\emptyset \Big)  &\le& \sum_{ k=n+1}^{2n} \e \Big[ \sum_{ \vert u \vert =k}  1_{ (s \le V(u) \le s+K, V(u_i) \ge a_i^{(n, \lambda)} , \forall i \le k)} \Big] \\
  	&=&  \sum_{ k=n+1}^{2n} \e \Big[ e^{ S_k} 1_{( s \le S_k \le s+K, S_i \ge a_i^{(n, \lambda)} , \forall i \le k )}\Big] \\
	&\le & \sum_{ k=n+1}^{2n} e^{ s+K}\, \p\Big( s \le S_k \le s+K, S_i \ge a_i^{(n, \lambda)}, \forall i \le k  \Big) . \end{eqnarray*}
 
 \noindent By \eqref{AS2}, $ \p\big( s \le S_k \le s+K, S_i \ge a_i^{(n, \lambda)} , \forall i \le k \big) \le c_{13}n^{-3/2}$ for all $n < k \le 2n$. Hence $ \p \big( \bigcup_{ k=n+1}^{2n}   E^{(n, \lambda)}_k  \not=\emptyset\big) \le c_{13}e^{- \lambda +K}$ proving the upper bound in \eqref{efk1}. $\Box$
 
\medskip

Using the notations in Lemma \ref{l:as} with the constant  $K$, we define for $n \ge2$ and $0\le \lambda \le {1\over 3} \log n$,   \begin{equation}\label{defanl}  A(n, \lambda):=  \Big\{ \cup_{ k=n+1}^{2n} \big( E^{(n, \lambda)}_k \cap F^{(n, \lambda)}_k \big) \not=\emptyset\Big\}. \end{equation} 
 
\noindent The following estimate will be useful  in the application  of Borel-Cantelli's lemma:

\begin{lemma}\label{l:bc}  There exists some constant $c_{14}>0$ such that for any $ n \ge2, 0\le \lambda \le {1\over 3}\log n$ and $m \ge 4 n, 0\le \mu \le {1\over 3} \log m$,  $$ \p \Big( A(n, \lambda) \cap A(m, \mu)\Big) \le c_{14}\, e^{ -\lambda - \mu}+ c _{14}\,e^{-\mu} { \log n \over \sqrt n}.$$
\end{lemma}

{ \noindent \bf Proof of Lemma \ref{l:bc}.} As   we mentioned before, the arguments that we use are very close to the computation of the second moment in the proof of Lemma 6.3 in \cite{AS12}.  The introduction of the events $F_k^{(n, \lambda)}$ in $A(n, \lambda)$, sometimes called a truncation argument, is   necessary to control the second moment:  the event  $F_k^{(n, \lambda)}$ keeps   the path $(V(u_i),  0 \le i \le k)$ of a particle $u $ in $E_k^{(n, \lambda)}$ to stay far away from  $(a_i^{(n, \lambda)}, 0 \le i \le k)$, otherwise the particle  $u$ would give a too large expectation in the second moment. Such truncation argument was already introduced in A\"{i}d\'{e}kon \cite{A11}.

 Let us enter into the details of the proof of Lemma \ref{l:bc}. Write for brevity   $$s:= { 1\over 2} \log n - \lambda, \qquad t:= { 1\over 2} \log m - \mu.$$  

\noindent
Similarly to \eqref{bk1} and \eqref{bk2}, we may construct a new probability $\q$ such that for all $n\ge1$,  ${ d \q \over d \p }\big\vert_{ \F_n}= W_n$,  $ \q \big( \xi_n = u \big \vert \F_n)={  e^{-V(u)} \over W_n}, \forall \vert u \vert =n$. Moreover under $\q$, $(V(\xi_n), n\ge 0)$ is distributed as the random walk $(S_n, n\ge0)$ defined in Section \ref{s:2},  and the spine decomposition similar to  (iii) in Section \ref{s:2} holds under $\q$. We refer to  \cite{BK04,CRW,LPP, AS12, Shi} for details.  It follows that \begin{eqnarray}
\p \Big( A(n, \lambda) \cap A(m, \mu)\Big) & \le& \e \Big[1_{ A(n, \lambda)} \, \sum_{ k=m+1}^{ 2m} \sum_{ \vert u \vert =k} 1_{( u \in E_k^{(m, \mu)} \cap F_k^{(m, \mu)})}\Big] \nonumber \\
	&=&\sum_{ k=m+1}^{ 2m}  \e_{\q}\Big[ 1_{ A(n, \lambda)} e^{ V(\xi_k)} 1_{( \xi_k   \in E_k^{(m, \mu)} \cap F_k^{(m, \mu)})}\Big] \nonumber \\
	&\le& e^{ t+K} \sum_{ k=m+1}^{ 2m}  \e_{\q} \Big[   A(n, \lambda) ,\,   \xi_k   \in E_k^{(m, \mu)} \cap F_k^{(m, \mu)} \Big] \nonumber \\
	&\le&  e^{ t+K} \sum_{ k=m+1}^{ 2m} \sum_{ l= n+1}^{2n}   \e_{\q}\Big[  \sum_{ \vert v\vert = l} 1_{(  v \in   E_l^{(n, \lambda)} \cap F_l^{(n, \lambda)} , \,   \xi_k   \in E_k^{(m, \mu)} \cap F_k^{(m, \mu)})} \Big] \nonumber\\
	&=:&    e^{ t+K} \sum_{ k=m+1}^{ 2m} \sum_{ l= n+1}^{2n}   I_{\eqref{anam}}(k, l)  . \label{anam}
\end{eqnarray}

\noindent For $n<l\le 2n \le {m \over 2}< k \le 2m$, we may decompose the   sum  on $ \vert v \vert = l$ as follows: $$\sum_{ \vert v \vert = l} 1_{(  v \in   E_l^{(n, \lambda)} \cap F_l^{(n, \lambda)} )}  = 1_{( \xi_l \in   E_l^{(n, \lambda)} \cap F_l^{(n, \lambda)}  )} +  \sum_{ p=1}^l \sum_{ u \in \Upsilon(\xi_p)} \sum_{ v \in \T(u), \vert v\vert_u= l-p}  1_{(  v \in   E_l^{(n, \lambda)} \cap F_l^{(n, \lambda)} )},$$

\noindent where $\T(u)$ denotes the subtree of $\T$ rooted at $u$ and $\vert v \vert_u= \vert v \vert - \vert u \vert $ the relative generation 
of $v \in \T(u)$.   Then  \begin{eqnarray}
&&   I_{\eqref{anam}}(k, l)  \nonumber \\
&=& \q\Big( \xi_l \in   E_l^{(n, \lambda)} \cap F_l^{(n, \lambda)},   \xi_k   \in E_k^{(m, \mu)} \cap F_k^{(m, \mu)} \Big)   + \sum_{ p=1}^l  \e_{\q} \Big[1_{( \xi_k   \in E_k^{(m, \mu)} \cap F_k^{(m, \mu)} )} \sum_{ u \in \Upsilon(\xi_p)} f_{k, l, p}(V(u)) \Big]  \nonumber\\
&=:& I_{\eqref{klp}}(k, l) +   \sum_{p=1}^l J_{\eqref{klp}}(k, l,p) ,  \label{klp}
\end{eqnarray} 

\noindent with $$  f_{k, l, p}(x)  :=  \e_{\q}  \Big[  \sum_{ v \in \T(u), \vert v\vert_u= l-p}  1_{(  v \in   E_l^{(n, \lambda)} \cap F_l^{(n, \lambda)} )} \big \vert V(u)=x\Big], \qquad x \in \r.$$

In what follows, we shall at first estimate $J_{\eqref{klp}}(k, l,p)$  then $I_{\eqref{klp}}(k, l)$. 
By the branching property at $u$ and by removing   the event $F_l^{(n, \lambda)} $ from the indicator function in $f_{k, l, p}(r)$, we get that 
 \begin{eqnarray} f_{k, l, p}(x)  
 &\le &   \e_x \Big[ \sum_{ \vert v \vert = l-p} 1_{( s \le V(v) \le s+K, V(v_i) \ge a_{i+p}^{(n,\lambda)} , \forall \, 0\le i \le l-p)}\Big] \nonumber \\
 &=&  e^{-x} \,  \e_x \Big[ e^{S_{l-p}}  1_{( s \le S_{l-p}  \le s+K,  S_i \ge a_{i+p}^{(n,\lambda)} , \forall \, 0\le i \le l-p)}\Big] \nonumber \\
 &\le & e^{-x + s+K} \,  \p_x \Big(  s \le S_{l-p}  \le s+K,  S_i \ge a_{i+p}^{(n,\lambda)} , \forall \, 0\le i \le l-p \Big), \label{fklp}
 \end{eqnarray}

\noindent where to get  the above equality, we applied an obvious modification of \eqref{many} for $\e_x$ instead of $\e$.

Let us denote by $\eqref{fklp}_{k, l,p}$ the probability term in \eqref{fklp}. To estimate $\eqref{fklp}_{k, l, p}$, 
we distinguish as in \cite{AS12} two cases: $p \le {n\over 2}$ and $ { n\over2} < p \le l$. Recall that $n<l\le 2n \le {m \over 2}< k \le 2m$.  If $p \le { n\over 2}$, 
$$ \eqref{fklp}_{k, l, p} \le  1_{ (x \ge0)} c_{15}\, { 1+x \over (l-p )^{ 3/2}}, $$ by using \eqref{AS2}.  Then for $1 \le p \le {n\over2}$,  $$ f_{k, l,p} (x) \le c_{15} 1_{(x\ge0)} e^{s+K-x} (1+x)  (l-p)^{-3/2}.$$

\noindent 
It follows that for all $n < l \le 2n , m < k \le 2m$,  \begin{eqnarray*}
\sum_{1 \le p \le n/2} J_{\eqref{klp}}(k, l,p) &\le & \sum_{ p=1}^{n/2}  \e_{\q} \Big[1_{( \xi_k   \in E_k^{(m, \mu)} \cap F_k^{(m, \mu)} )} \sum_{ u \in \Upsilon(\xi_p)}  c_{15} 1_{( V(u) \ge0)} e^{s+K-V(u)} { 1+V(u) \over (l-p)^{3/2}}  \Big]   \\
	&\le& c_{16}\,  e^{s} n^{-3/2}\, \sum_{ p=1}^{n/2} \e_{\q} \Big[1_{( \xi_k   \in E_k^{(m, \mu)} \cap F_k^{(m, \mu)} )} \sum_{ u \in \Upsilon(\xi_p)}   1_{( V(u) \ge0)} e^{ -V(u)} ( 1+V(u) ) \Big] \\
	&\le& c_{16}\, K\,  e^{s} n^{-3/2}\, \sum_{ p=1}^{n/2} \e_{\q}\Big[1_{( \xi_k   \in E_k^{(m, \mu)} \cap F_k^{(m, \mu)} )} e^{ - (p-1)^{1/12}}\Big] , \end{eqnarray*}

\noindent where the last inequality is due to the definition of  $\xi_k \in F_k^{(m, \mu)}$ [noticing that $a_p^{(m, \mu)}=0$ and $b_p^{(k, m)}= p^{1/12}$ for all $p \le n/2< m/2$].   Then we get that  \begin{equation}\label{I3a} \sum_{1 \le p \le n/2} J_{\eqref{klp}}(k, l,p)  \le c_{17} e^{s} n^{-3/2}\, \q \Big( \xi_k   \in E_k^{(m, \mu)} \Big) \le c_{18} e^{s} n^{-3/2}\,  m^{-3/2},
\end{equation}

\noindent since $\q \Big( \xi_k   \in E_k^{(m, \mu)} \Big)= \p( t \le S_k \le t+K, S_i \ge a_i^{(m, \mu)}, \forall \, 0\le i \le k\Big) \le c_{19} m^{-3/2}$ for all $m< k \le 2m$, by using \eqref{AS2}.

 Now considering $ { n \over2} < p \le l$, $a^{(n,\lambda)}_{i+p}= s$ for any $0 \le i \le l-p$, hence $$ \eqref{fklp}_{k, l, p} = 1_{(x \ge s)} \p_x \Big( s \le S_{l-p} \le s+K, \underline S_{l-p} \ge s\Big) \le 1_{(x\ge s)} c_2 (1+K)^2\, { 1+x-s\over (1+l-p)^{3/2}},$$ by \eqref{AS1}.  It follows   that \begin{eqnarray*}
  \sum_{ {n\over2} < p \le l} J_{\eqref{klp}}(k, l,p) & \le&   \sum_{ {n\over2} < p \le l}  \e_{\q}\Big[1_{( \xi_k   \in E_k^{(m, \mu)} \cap F_k^{(m, \mu)} )} \sum_{ u \in \Upsilon(\xi_p)}  c_2 (1+K)^21_{( V(u) \ge s)} e^{s+K-V(u)} { 1+V(u)-s \over (1+l-p)^{3/2}}  \Big]  . 
\end{eqnarray*}

\noindent 
By the definition of    $\xi_k \in F_k^{(m, \mu)}$,  for all $p \le l \le 2n \le {m\over2}$, we have that  $$\sum_{ u \in \Upsilon(\xi_p)}    1_{( V(u) \ge s)} e^{ -V(u)} (1+V(u)-s) \le \sum_{ u \in \Upsilon(\xi_p)}    1_{( V(u) \ge 0)} e^{ -V(u)} (1+V(u) ) \le K e^{- (p-1)^{1/12}}  .$$

\noindent 
Then \begin{equation}\label{I3b}  \sum_{   {n\over2} \le p \le l} J_{\eqref{klp}}(k, l,p)  \le c_{20}e^s\,  e^{- n^{1/13}} \q(\xi_k   \in E_k^{(m, \mu)}) \le c_{21}  e^{s - n^{1/13}} m^{-3/2} . \end{equation}

 \noindent Combining \eqref{I3a} and \eqref{I3b}, we get that \begin{equation}\label{I3ab}  \sum_{   1 \le p \le l} J_{\eqref{klp}}(k, l,p) \le (c_{18} n^{-3/2}+ c_{21} e^{- n^{1/13}}) e^s \, m^{-3/2}. \end{equation}
 
 It remains to estimate $I_{\eqref{klp}}(k, l)$ for $n <l \le 2n$ and $m< k \le 2m$ [in particular $l < k$]. We have \begin{eqnarray}
&&   I_{\eqref{klp}}(k, l) \nonumber \\ &\le & \q\Big( \xi_l \in   E_l^{(n, \lambda)} ,   \xi_k   \in E_k^{(m, \mu)}   \Big)  \nonumber \\
	&    =&      \p \Big( s \le S_l \le s+K, S_i \ge a_i^{(n,\lambda)} , \forall \, i \le l, \, t \le S_k \le t+K, S_j \ge a_j^{(m, \mu)}, \forall \,   j \le k\Big)  . \label{I2a}
\end{eqnarray}

\noindent Let us denote by $ \eqref{I2a}_{k, l}$ the probability term in \eqref{I2a}.
Using the Markov property at $l$, we get that  \begin{eqnarray*}
\eqref{I2a}_{k, l}   &=& \e \Big[ 1_{( s \le S_l \le s+K, S_i \ge a_i^{(n,\lambda)} , \forall \, 0 \le i \le l)} \p_{S_l} \Big( t \le S_{k-l} \le t+K, S_j \ge a_{j+l}^{(m, \mu)}, \forall \, 0\le j \le k-l\Big) \Big]  \\
	&\le& {c_{21}  \over (k-l)^{3/2}} \e \Big[ 1_{( s \le S_l \le s+K, S_i \ge a_i^{(n,\lambda)} , \forall \, 0 \le i \le l)}  (1+S_l) \Big] \qquad \mbox{(by \eqref{AS1})} \\
	&\le& c_{22} (1+s+K)  (k-l)^{-3/2} \, l^{-3/2}.
\end{eqnarray*}

\noindent Based on the above estimate and \eqref{I3ab}, we deduce from  \eqref{anam} and \eqref{klp}  that \begin{eqnarray*}
&&\p \Big( A(n, \lambda) \cap A(m, \mu)\Big) \\ & \le& c_{23} \, e^{ t+K} \sum_{ k=m+1}^{ 2m} \sum_{ l= n+1}^{2n}  \Big(  (1+s+K)  (k-l)^{-3/2} \, l^{-3/2} + e^s\,  e^{- n^{1/13}} m^{-3/2}+ e^{s} n^{-3/2}\,  m^{-3/2}\Big) \\
	& \le& c_{24}  e^{-\lambda -\mu} + c_{24} e^{-\mu} { \log n \over \sqrt n},
\end{eqnarray*} proving the Lemma. $\Box$

 \medskip
 
\noindent{\bf Proof of the divergence  part in Theorem \ref{t:min}.} Let $f$ be nondecreasing such that $\int^\infty {dt \over t e^{f(t)}} =\infty$. Without any loss of generality we may assume that $ \sqrt{ \log t} \le  e^{f(t)}  \le ( \log t)^2$ for all large $t \ge t_0$ (see e.g.  \cite{erdos} for a similar justification).  
Denote by $$B_x(k):=\Big\{ \M_n +x \le { 1 \over 2} \log n -   f(n+k), \, \, \mbox{i.o. as }  n \to \infty \Big\}, \qquad x\in \r, \, k \ge0.$$

Let us first prove that there exists some constant  $c_{25}>0$ such that for any $x \in \r$ and $k \ge 0$, \begin{equation} \p \Big(  B_x(k)\Big) \ge c_{25}. \label{bc2}
\end{equation}

To this end, we take $n_i:= 2^i $ for $i \ge1$, $\lambda_i:= f(n_{i+1} +k) +x  +K $ and consider the event $A_i:=A(n_i, \lambda_i)$ in \eqref{defanl}.  There is some integer $i_0\equiv i_0(x, k)\ge1$ such that for all $i \ge i_0$, $0\le \lambda_i \le {1\over 3} \log n_i$.   By Lemma \ref{l:as},  $$ c_{12} \, e^{ - \lambda_i} \le \p( A_i) \le {1\over c_{12}} \, e^{- \lambda_i}, \qquad \forall i \ge i_0. $$ 

\noindent Note that  $\int^\infty { dt \over t e^{f(t +k)}} \ge \int^\infty { dt \over (t+k) e^{f(t+k)}}=\infty$, and   $\int_{n_{i+1}}^{n_{i+2} }  { dt \over t e^{f(t +k)}}  \le ( \log 2) e^{- f(n_{i+1} +k)}$ by the monotonicity  of $f$.  Hence $\sum_i e^{-\lambda_i}=\infty$. By Lemma \ref{l:bc}, we have for any $i \ge i_0$ and $j \ge i+2$, $$ \p\Big( A_i \cap A_j\Big) \le c_{14}\, e^{ - \lambda_i -\lambda_j} + c_{14} e^{- \lambda_j } { \log n_i \over \sqrt{n_i}},$$ 

\noindent which implies that $ \sum_{ i, j=i_0}^k  \p\Big( A_i \cap A_j\Big) \le c_{14} \big( \sum_{i=i_0}^k e^{ - \lambda_i}\big)^2 + 2 c_{14} \big( \sum_{i=i_0}^k e^{ - \lambda_i}\big) \times \big( \sum_{i=1}^\infty { \log n_i \over \sqrt{ n_i}}\big). $ Using the lower bound $\p(A_i) \ge c_{12} \, e^{ - \lambda_i}$ and the fact that    $\sum_{i=i_0}^k e^{-\lambda_i} \to \infty $ as $k \to \infty$, we  obtain   that $$ \limsup_{k\to \infty} { \sum_{1 \le  i , j \le k} \p(A_i \cap A_j) \over \big[ \sum_{ i=1}^k \p(A_i)\big]^2} \le {c_{14} \over c_{12}^2}.$$

 By Kochen and Stone \cite{KS}'s version of the Borel-Cantelli lemma,     $\p( A_i, \mbox{i.o.} \, i \to \infty) \ge c_{12}^2 / c_{14}=: c_{25}$ which does not depend  on $(x, k)$.   Observe that $\{ A_i, \mbox{i.o.} \, i \to \infty\} \subset B_x(k)$, in fact, for those $i $ such that $A_i\equiv A(n_i, \lambda_i)$ holds, by the definition   \eqref{defanl}, there exits some $n \in ( n_i, n_{i+1}] $ such that $\M_n \le { 1 \over 2} \log n_i - \lambda_i + K = { 1 \over 2} \log n_i - f(n_{i+1} +k)  -x \le  { 1 \over 2} \log n - f(n  +k)  -x$.   Hence we get \eqref{bc2}.

  We have proved  that for any $x\in \r$ and $k\ge0$, $\p(B_x(k)) \ge c_{25}$. For any $k\ge0$, the events $B_x(k)$ are non-increasing on $x$. Let $B_\infty (k):= \cap_{ i=1}^\infty B_i (k) $ [then $B_\infty(k)$ is  nothing but $\{ \liminf_{n \to \infty} ( \M_n -{1\over2} \log n + f(n+k))= - \infty\}$]. By the monotone convergence, $ \p(B_\infty(k)) \ge c_{25}$, for all $k\ge0$.  Moreover, for any $x\in \r$, $ \p_x(B_\infty(k))= \p(B_\infty(k))\ge c_{25}$. On the other hand, if we denote by  $Z_k:= \sum_{\vert u \vert =k} 1$ the number of particles in the $k$-th generation, then by the branching property,  $$ \p \Big( B_\infty (0)\, \big\vert \, \F_k \Big) =1_{( Z_k >0)} \Big( 1- \prod_{ \vert u \vert =k} ( 1- \p_{V(u)}(B_\infty(k))) \Big) \ge 1_{(Z_k >0)} \Big( 1- (1-c_{25})^{Z_k}\Big).$$

It is well-known (cf. \cite{AN}, pp.8) that     $\S = \{ \lim_{k \to \infty} Z_k = \infty\}$.  Then  by letting $k \to \infty$ in the above inequality, we get that  $$1_{B_\infty(0)} = \lim_{k \to \infty} \p \Big( B_\infty (0)\, \big\vert \, \F_k \Big) \ge 1_{\S}, \qquad \p\mbox{-a.s.}$$

\noindent Clearly $\S^c \subset B_\infty(0)^c$ by the convention on the definition of $\M_n$ on $\S^c$. Hence   $\S=B_\infty(0)$,  $\p$-a.s.   This proves the divergence part of Theorem \ref{t:min}. $\Box$

\section{Proof of Proposition \ref{p:lower}}\label{s:4}

The main technical part was already done in A{\"{\i}}d{\'e}kon and Shi \cite{AS12}:  
\begin{lemma} [\cite{AS12}]  \label{l:AS41} Assume  \eqref{hyp}  and \eqref{int3}.  For any fixed $\alpha\ge0$, there exist some $\delta= \delta(\varepsilon_0) >0$ and $c_{26}=c_{26}(\alpha, \delta)>0$ such that for all $n \ge2$,  \begin{equation} \mbox{Var}_{ \qa} \left( { \sqrt n W_n^{(\alpha)}\over D_n^{(\alpha)}}\right) \le c_{26}\,\left(  n^{-\delta} + \sup_{ k_n^{1/3} \le x \le k_n } \Big\vert  { h_{x+\alpha} ( n-k_n) \over h_\alpha(n) }- 1 \Big\vert  \right),
\end{equation}  where $k_n:= \lfloor n^{1/3}\rfloor $ and  $  h_x(j):= { \sqrt j \, \p_x ( \underline S_j \ge 0 ) \over R(x)}$ for $  j\ge1, \, x\ge0.$ 
\end{lemma}

{\noindent \bf Proof of Lemma \ref{l:AS41}.}  The Lemma   was implicitly contained   in the proof of Proposition 4.1 in \cite{AS12}. In fact,  in their proof of the convergence that $\mbox{Var}_{ \qa} \left( { \sqrt n W_n^{(\alpha)}\over D_n^{(\alpha)}}\right) \to0$,  we choose  $k_n:= \lfloor n^{1/3}\rfloor$ in their definition of $E_n:= E_{n,1} \cap E_{n, 2} \cap E_{n,3}$ (see the equation (4.6) in \cite{AS12}, Section 4).   We claim that for some constant $\delta_1=\delta_1(\varepsilon_0)>0$,  there is some $c_{27}=c_{27}(\delta_1, \alpha)>0$ such that for all $n \ge1$, \begin{eqnarray}\label{en123} \qa \Big( E_n^c\Big) &\le  &c_{27}\, n^{-\delta_1}, \\
    \sup_{ k_n^{1/3} \le x \le k_n} \qa \Big( E^c_n \big \vert V(\xi_{k_n} ) =x\Big) & \le& c_{27} n^{-\delta_1}. \label{en123b}
\end{eqnarray}

In fact, according to the definition of $E_{n,1}$ in \cite{AS12}, $$ \qa \Big( E_{n, 1}^c\Big) \le \qa \Big( \{V(\xi_{k_n}) > k_n \} \cup\{ V(\xi_{k_n}) < k_n^{1/3}\}\Big)  + \sup_{ k_n^{1/3} \le  x \le k_n} \qa_x\Big( \cup_{i=0}^{n-k_n}  \{ V(\xi_i) < k_n^{1/6} \} \Big).$$

\noindent By \eqref{bk4}, $$ \qa \Big( V(\xi_{k_n}) < k_n^{1/3}\Big) = { 1\over R_\alpha(0)} \e \Big( 1_{( \underline S_{k_n} \ge -\alpha, S_{k_n} < k_n^{1/3})} R_\alpha(S_{k_n})\Big)  \le {R_\alpha(k_n^{1/3})\over R_\alpha(0)}\, \p( \underline S_{k_n} \ge -\alpha) \le c_{28} k_n^{-1/6},$$ 

\noindent and $$ \qa \Big( V(\xi_{k_n}) > k_n \Big) \le \e \Big[ 1_{ ( S_{k_n} > k_n)} R_\alpha(S_{k_n})\Big] \le   \sqrt{ \p\Big( S_{k_n} > k_n \Big)}  \sqrt{ \e[  R_\alpha(S_{k_n})^2]} \le c_{29} k_n \, \sqrt{ \p\Big( S_{k_n} > k_n \Big)} ,$$

\noindent since $R_\alpha(x) \sim c_R  x $ as $x \to \infty$.   The condition \eqref{int3} ensures that $\e  ( \vert S_1\vert^{2+\varepsilon_0})< \infty$ which in turn  implies that $\e ( \vert S_k\vert ^{ 2+\varepsilon_0}) \le c_{30} k^{ 1+ \varepsilon_0/2}$ for any $k \ge1$ (see Petrov \cite{P95}, pp.60). Hence  $$ \qa \Big( V(\xi_{k_n}) > k_n \Big) \le c_{31} k_n ^{ - \varepsilon_0/4}.$$

Now for $k_n^{1/3} \le x \le k_n$,  let $\tau= \inf\{ i \ge0: S_i < k_n^{1/6}\}$, then  the absolute  continuity \eqref{bk4} at $\tau$ reads as  \begin{eqnarray*}
 \qa_x\Big( \cup_{i=0}^{n-k_n}  \{ V(\xi_i) < k_n^{1/6} \} \Big) &=& { 1\over R_\alpha(x) } \e_x \Big[ 1_{ (  \tau \le n- k_n)} R_\alpha(S_\tau) 1_{( \underline S_ \tau\ge -\alpha)}\Big] \\
 	&\le &  { R_\alpha(k_n^{1/6}) \over R_\alpha(x) }   \le   { R_\alpha(k_n^{1/6}) \over R_\alpha(k_n^{1/3}) } \le    c_{32}\, k_n^{ -1/6},
\end{eqnarray*} since $x\ge k_n^{1/3}$.  Assembling the above estimates yields that $$ \qa \Big( E_{n, 1}^c\Big) \le c_{33} k_n^{-\varepsilon_0/4},$$

\noindent [we may assume $\varepsilon_0 \le 2/3$].  Let us follow the proof of Lemma 4.7 in \cite{AS12},   we remark that on $E_{n, 1}$, $ V(\xi_i) \ge k_n^{1/6}$ for all $k_n \le i \le n$,   and it was shown in \cite{AS12}   that $$ \qa_x\Big( E_{n, 1} \cap E_{n, 2}^c \Big) \le \sum_{ i=k_n}^{n-1}  \e _x \Big[    1_{( \eta + \widetilde \eta > e^{S_i /2})}  \big[ \eta + { \widetilde \eta \over S_i+\alpha+1} \big] \, 1_{(S_i \ge  k_n^{1/6})}\Big].$$

 By the integrability assumption \eqref{int3}, since $\widetilde \eta= \int_0^\infty x e^{-x} \L(dx) \le \sqrt{ \eta \, \int_0^\infty x^2 e^{-x} \L(dx)}$, it is easy to see that $ \e( \widetilde \eta^p)< \infty$ for some $p >1$. It follows that 
 $  \qa_x\Big( E_{n, 1} \cap E_{n, 2}^c \Big) \le n \, \e    \Big(    1_{( \eta + \widetilde \eta > e^{k_n^{1/6} /2})}  \big[ \eta +   \widetilde \eta \big]  \Big) \le  c_{34} n^{-10}$ .    Finally by (4.9) in \cite{AS12}, $\qa( E_{n,1}\cap E_{n,2} \cap E^c_{n, 3}) \le c_{35} n^{-10}$, hence we get \eqref{en123}. 

From  \eqref{en123}, it suffices to follow line-by-line  the proof of Proposition 4.1 in \cite{AS12}: In Lemma 4.4 of \cite{AS12}, we can get $n^{ - 1- \delta_1/4}$ instead of $o({1\over n})$ [by replacing in its proof $\varepsilon$ by $n^{-\delta_1/4}$]. In their proof of Lemma 4.5, taking $\eta_1= {1\over n}$ and we arrive at $$ \e_{\qa}\left[\Big( \sqrt n { W_n^{(\alpha)} \over D_n^{(\alpha)}}\Big)\right]^2 \le c_{36} n^{-\delta_1/4} + (1+O({1\over n})) \e_{\qa}\big[ \sqrt n { W_n^{(\alpha)} \over D_n^{(\alpha)}}\big] \, \sup_{ k_n^{1/3} \le x \le k_n} { \sqrt n \p( \underline S_{n-k_n} \ge - \alpha- x ) \over R_\alpha(x)}.$$

 \noindent The Lemma follows    because  $  \e_{\qa}\big[ \sqrt n { W_n^{(\alpha)} \over D_n^{(\alpha)}}\big]= h_\alpha(n) $, and $h_\alpha(n) \to \theta$ when $n \to \infty$,  as shown in \cite{AS12}.   
$\Box$

\medskip

{\noindent\bf Proof of Proposition \ref{p:lower}.}   It is enough to  prove that for any $\alpha \ge0$,  \begin{equation}
\liminf_{n \to \infty} \sqrt n { W_n^{(\alpha)} \over D_n^{(\alpha)}} = \theta, \qquad \mbox{$\qa$-a.s.}, \label{last1}
\end{equation} 

\noindent where $\theta$ is defined in \eqref{K1}. In fact,  under \eqref{hyp}, \eqref{int1} and  \eqref{int2}, $D_n^{(\alpha)}$ converges in mean to $D_\infty^{(\alpha)}$ (see  \cite{Shi}, Chapter 5, also see \cite{BK04}, Theorem 10.2 (i) with an extra $\log\log\log$-term). Then on $\{D_\infty^{(\alpha)}>0\}$, $\p$ and $\qa$ are equivalent. Moreover, 
 as shown in \cite{AS12},  $\p$-almost surely on $\{ \inf_{ \vert u \vert \ge0} V( u) \ge - \alpha\}$,  $W_n^{(\alpha)}=W_n$ and $ \lim_{n \to \infty} D_n^{(\alpha)} =  c_R \, D_ \infty$, therefore Proposition \ref{p:lower} follows easily  from \eqref{last1}.

Now to prove \eqref{last1}, since $\sqrt n { W_n^{(\alpha)} \over D_n^{(\alpha)}}  \to \theta$ in probability under $\qa$  (\cite{AS12}), it suffices  to prove that \begin{equation}
\liminf_{n \to \infty} \sqrt n { W_n^{(\alpha)} \over D_n^{(\alpha)}} \ge  \theta, \qquad \mbox{$\qa$-a.s.} \label{last2}
\end{equation}

\noindent To this end,  using  Lemmas \ref{l:AS41} and \ref{l:rw} we get  some constant $\delta_2>0$ such that for all $n \ge n_0$, \begin{equation}
\mbox{Var}_{\qa} \Big( \sqrt n { W_n^{(\alpha)} \over D_n^{(\alpha)}}  \Big) \le n^{-\delta_2}. \label{last3}
\end{equation}

\noindent Let $n_j:= j^{-3/\delta_2}$ for $j \ge j_0$ and choose an arbitrary small  $\varepsilon>0$. We are going to show that \begin{equation}
\sum_{j\ge j_0} \,   \qa  \Big( \inf_{ n_j \le n \le n_{j+1}} \sqrt n { W_n^{(\alpha)} \over D_n^{(\alpha)}}  < (1-\varepsilon) \theta\Big) < \infty, \label{last4}
\end{equation} from which the Borel-Cantelli lemma   yields  \eqref{last2}. 

To prove \eqref{last4}, let $\widehat \F_n:= \F_n \vee \G_n$, where $\G_n$, defined in Section \ref{s:2},  denotes the $\sigma$-fields generated by the spine up to generation $n$.  Then $\qa$-a.s.,   \begin{eqnarray*}
\e_{\qa} \left [ {1\over R_\alpha(V(\xi_{n_{j+1}}))} \Big\vert\, \widehat F_n\right] &= &\e_{\qa_{V(\xi_n)}} \Big [ {1\over R_\alpha( V(\xi_{n_{j+1}- n}))}\Big] \\
	&=& {1\over R_\alpha( V(\xi_n))}  \p_{V(\xi_n)} \Big( \underline S_{ n_{j+1}-n}  \ge -\alpha\Big)  \qquad (\mbox{by }\eqref{bk4})\\
	&\le&  {1\over R_\alpha( V(\xi_n))} .
\end{eqnarray*}

\noindent It follows that for all $n \le n_{j+1}$, $$ \e _{ \qa} \left [  {1\over R_\alpha(V(\xi_{n_{j+1}}))} \Big\vert\,  \F_n\right ]  \le \e _{ \qa} \left [ {1\over R_\alpha(V(\xi_n))} \Big\vert\,  \F_n\right ]   = {W_n^{(\alpha)} \over D_n^{(\alpha)}},$$ where the last equality comes from Lemma 4.2 in \cite{AS12}. 
 Consequently for all $n_j \le n \le n_{j+1}$, $$ \sqrt n \,  {W_n^{(\alpha)} \over D_n^{(\alpha)}} \ge Y_n:= \sqrt{n_j} \, \e _{ \qa} \left [   {1\over R_\alpha(V(\xi_{n_{j+1}}))} \Big\vert\,  \F_n\right ]  .$$

Remark that $(Y_n, n_j \le  n \le n_{j+1})$ is a martingale with mean $\e_{\qa}(Y_{n_j})= \sqrt{ n_j} \, \e_{\qa} ({ 1\over R_\alpha(S_{n_j})}) \ge (1-\varepsilon) \theta$. The Doob $L^2$-inequality implies that \begin{eqnarray*}
\qa \Big( \max_{ n_j \le n \le n_{j+1}} \big\vert Y_n - \e_{\qa}(Y_{n_j})\big\vert \ge{ \varepsilon \over 2}\theta \Big)  &\le &{4\over \varepsilon^2\theta^2} \mbox{Var}_{\qa} ( Y_{n_{j+1}})  \\ 
&\le& c_{36}\, n_j^{-\delta_2}  = c_{36}  j^{-3},
\end{eqnarray*} by \eqref{last3} and the fact that $Y_{n_{j+1}}= \sqrt{n_j \over n_{j+1}}\, {W_{n_{j+1}}^{(\alpha)} \over D_{n_{j+1}}^{(\alpha)}}$.  Finally for all large $j$, $$\qa  \Big( \inf_{ n_j \le n \le n_{j+1}} \sqrt n { W_n^{(\alpha)} \over D_n^{(\alpha)}}  < (1-\varepsilon) \theta\Big)  \le 
\qa \Big( \max_{ n_j \le n \le n_{j+1}} \big\vert Y_n - \e_{\qa}(Y_{n_j})\big\vert  \ge { \varepsilon \over 2}\theta \Big) \le c_{36} \, j^{-3} ,$$

\noindent proving \eqref{last4} and then completing the proof of Proposition \ref{p:lower}. $\Box$

\bigskip
\noindent{\bf Acknowledgements.}  I am grateful to an anonymous referee for her/his helpful suggestions  and careful reading  of  the manuscript. I also thank Vladimir Vatutin for sending me the paper Eppel \cite{E79}.


\begin{thebibliography}{12}
\baselineskip=12pt

 
 
 \bibitem{AbR}   Addario-Berry, L. and Reed, B. (2009).
 Minima in branching random walks.
  {\it Ann.  Probab.}, {\bf 37},  {1044--1079}.


\bibitem{A11} A{\"{\i}}d{\'e}kon, E. (2011+).
     Convergence in law of the minimum of a branching random walk.      {\it Ann. Probab., to appear.}
 
 \bibitem{ABBS} A{\"{\i}}d{\'e}kon, E., Berestycki, J., Brunet, E. and Shi, Z. (2011+) The branching Brownian motion seen from its tip. (preprint, http://arxiv.org/abs/1104.3738). 
 
 \bibitem{AJ11} A{\"{\i}}d{\'e}kon, E. and Jaffuel, B. (2011).  Survival of branching random walks with absorption. {\it Stochastic Process. Appl. }  {\bf 121}   pp. 1901--1937. 
 
\bibitem{AS12} A{\"{\i}}d{\'e}kon, E. and Shi, Z. (2012+). The Seneta-Heyde scaling for the branching random walk. {\it Ann. Probab., to appear}.


\bibitem{ABK}  Arguin, L.P., Bovier, A. and Kistler, N. (2011+) The Extremal Process of Branching Brownian Motion. (preprint, http://arxiv.org/abs/1103.2322).

 
 \bibitem{AN} Athreya, K.B and Ney, P.E.  (1972). {\it Branching processes.} Springer-Verlag, Berlin, New-York.  
 



\bibitem{B76}
    Biggins, J.D. (1976).
    The first- and last-birth problems for a
    multitype age-dependent branching process.
    {\it Adv. Appl. Probab.} {\bf 8}, 446--459.



\bibitem{BK04}
    Biggins, J.D. and Kyprianou, A.E. (2004).
    Measure change in multitype branching.
    {\it Adv. Appl. Probab.} {\bf 36}, 544--581.

\bibitem{BK05}
    Biggins, J.D. and Kyprianou, A.E. (2005).
    Fixed points of the smoothing transform: the
    boundary case.
    {\it Electron. J. Probab.} {\bf 10}, Paper
    no.~17, 609--631.

 
\bibitem{B78}
    Bramson, M.D. (1978).
    Maximal displacement of branching Brownian
    motion.
    {\it Comm. Pure Appl. Math.} {\bf 31}, 531--581.

 \bibitem{BZ06}
    Bramson, M.D. and Zeitouni, O. (2009).
    Tightness for a family of recursion equations. {\it Ann. Probab.,} {\bf 37},          615--653
    

\bibitem{CRW}
    Chauvin, B., Rouault, A. and Wakolbinger, A.
    (1991).
    Growing conditioned trees.
    {\it Stoch. Proc. Appl.} {\bf 39}, 117--130.



 




\bibitem{erdos} Erd\H{o}s P. (1942). On the law of the iterated logarithm.  {\it Ann. of Math. } (2) 43,   419--436. 

\bibitem{E79} 
Eppel, M.S. (1979).  A local limit theorem for the first overshoot. {\it Siberian Math. J.} {\bf 20}  130--138.


 \bibitem{feller}
   Feller, W. (1971).
    {\it An Introduction to Probability Theory and
   its Applications.} Vol. II. Second edition.
     Wiley, New York.

\bibitem{H74}
    Hammersley, J.M. (1974).
    Postulates for subadditive processes.
    {\it Ann. Probab.} {\bf 2}, 652--680.



\bibitem{HS09}
    Hu, Y.\ and Shi, Z.  (2009).
    Minimal position and critical martingale
    convergence in branching random walks,
    and directed polymers on disordered
    trees.
    {\it Ann.\ Probab.}, {\bf 37}, 742--789.

\bibitem{J12}
    Jaffuel, B.\ (2012+).
    The critical barrier for the survival of
    the branching random walk with absorption.
  To appear in  {\it Annales de l'Institut Henri Poincar\'e.}  {\tt ArXiv math.PR/0911.2227}.
  
  
 \bibitem{K75}
    Kingman, J.F.C. (1975).
    The first birth problem for an age-dependent
    branching process.
    {\it Ann. Probab.} {\bf 3}, 790--801.

\bibitem{KS}  Kochen, S. and Stone, C. (1964). A note on the Borel-Cantelli lemma. { \it Illinois J. Math.} { \bf 8}, 248--251. 

\bibitem{K76}
    Kozlov, M.V. (1976).
    The asymptotic behavior of the probability of
    non-extinction of critical branching processes
    in a random environment.
    {\it Theory Probab. Appl.} {\bf 21}, 791--804.




\bibitem{L97}
    Lyons, R. (1997).
    A simple path to Biggins' martingale convergence
    for branching random walk.
    In: {\it Classical and Modern Branching
    Processes} (Eds.: K.B.~Athreya and P.~Jagers).
    {\it IMA Volumes in Mathematics and its
    Applications} {\bf 84}, 217--221. Springer, New
    York.

\bibitem{LPP}
    Lyons, R., Pemantle, R. and Peres, Y. (1995).
    Conceptual proofs of $L\log L$ criteria for mean
    behavior of branching processes.
    {\it Ann. Probab.} {\bf 23}, 1125--1138.

\bibitem{M10} Madaule, Th. (2010+).  Convergence in law for the branching random walk seen from its tip. (Preprint, http://arxiv.org/abs/1107.2543).


 
 \bibitem{P95}
     Petrov, V.V. (1995).
    {\it Limit Theorems of Probability Theory.}
    Clarendon Press, Oxford.

\bibitem{Shi} Shi, Z. (2012+). {\it Branching random walks.}  Saint-Flour's summer course. 

\bibitem{VW09} Vatutin,  V. and Wachtel, V. (2009). Local probabilities for random walks conditioned
to stay positive. {\it Probab. Th. Rel. Fields.} {\bf 143} {177--217.}
 
\end{thebibliography}
\end{document}